\documentclass[openany, a4paper, 12pt]{article}

\usepackage{amsfonts}

 \usepackage{mathrsfs,amsfonts,amsmath}
 \usepackage{color}

 \makeatletter
 \@addtoreset{equation}{section}
 \makeatother
 \newtheorem{theorem}{Theorem}[section]
 \newtheorem{definition}{Definition}[section]
 \newtheorem{hypothesis}{Hypothesis}[section]
 \newtheorem{lemma}{Lemma}[section]
 \newtheorem{proposition}{Proposition}[section]
 \newtheorem{corollary}{Corollary}[section]
 \newtheorem{remark}{Remark}[section]
 \newtheorem{example}{Example}[section]

 \def\beqlb{\begin{eqnarray}}\def\eeqlb{\end{eqnarray}}
 \def\beqnn{\begin{eqnarray*}}\def\eeqnn{\end{eqnarray*}}
 \def\beqn{\begin{equation}}
 \def\eeqn{\end{equation}}

 \def\mbb{\mathbb}

\def\dz{\delta}

\def\ez{\epsilon}
\def\we{\wedge}

\def\ra{\rangle}
\def\la{\langle}

\begin{document}

\vskip1.0cm

\centerline{\large\textbf{Strong uniqueness for a class of singular
SDEs}}

\smallskip

\centerline{\large\textbf{ for catalytic branching diffusions
}\footnote{Supported by NSFC (No. 10721091)}}

\bigskip

\centerline{By Hui He\footnote{ \textit{E-mail address:} {\tt
hehui@bnu.edu.cn }} }

\medskip

\centerline{ School of Mathematical Sciences, Beijing Normal
University}

\smallskip

\centerline{ Beijing 100875, People's Republic of China}

\bigskip

{\narrower{\narrower{\narrower{\narrower

\begin{abstract}
A new result for the strong uniqueness for catalytic branching
diffusions is established, which improves the work of Dawson, D.A.;
Fleischmann, K.; Xiong, J.[Strong uniqueness
  for cyclically symbiotic branching diffusions. \textit{Statist. Probab.
  Lett.} \textbf{73}, no. 3, 251--257 (2005)].
\end{abstract}

\smallskip


\smallskip


\smallskip

\noindent\textbf{Abbreviated Title:} Uniqueness for catalytic
branching diffusions

\par}\par}\par}\par}

\bigskip\bigskip

\section{Introduction}
Stochastic differential equation(SDE) is a very important tool in
the theory of diffusion processes. Many investigations were devoted
to the problems of existence, uniqueness, and properties of
solutions of SDEs.  The well-known  result of Yamada and Watanabe
says that if a solution of a SDE exists and the pathwise uniqueness
of solutions holds, then the SDE admits a unique strong solution;
see Ikeda and Watanabe (1989, p.163) and Revuz and Yor (1991,
p.341). Then the study of pathwise uniqueness is of great interest.
For a long time much has been known about uniqueness for
one-dimensional stochastic differential equations (SDEs) with
singular coefficients. The diffusion coefficient can be
non-Lipschitz and degenerate; the drift can be singular and involve
local time. See, e.g., Cherny and Engelbert (2005) for a survey.
Especially, some results on pathwise uniqueness (strong uniqueness)
for SDEs have been obtained for certain H\"{o}lder continuous
diffusion coefficients; see Revuz and Yor (1991, Chapter IX-3) and
Ikeda and Watanabe (1989, p.168). These results are sharp; see
Barlow (1982). However, there are much less results on the pathwise
uniqueness beyond the Lipschitz (or locally Lipschitz) conditions in
the higher-dimensional case. Recent work in this direction includes
the papers of Fang and Zhang (2005), Swart (2001, 2002) and
DeBlassie (2004).

In this work, we shall study the pathwise uniqueness for a class of
degenerate stochastic differential equations with non-Lipschitz
coefficients. Our interest is motivated by models of catalytic
branching networks that include \textit{mutually catalytic
branching} and \textit{cyclically catalytic branching} diffusions;
see Dawson and Fleischmann (2000) for a survey on these systems. For
models with mutually catalytic branching and cyclically catalytic
branching,  the branching rate of one type is allowed to depend on
the frequency of the other types. The intuition is that the presence
of different types affects the branching of  other types. By the
interaction over all species, the basic independence assumption in
classical branching theory is violated. Uniqueness for those models
is usually hard to prove.  Recently, Athreya et al (2002), Bass and
Perkins (2003) and Dawson and Perkins (2006) studied \textit{weak
uniqueness} for
 \beqn
 \label{0.1}
 dX_t^i=b^i(X_t)dt+\sqrt{2\sigma_i(X_t)X_t^i}dB_t^i,
 ~~~i=1,2,\cdots,d
 \eeqn
in  $\mathbb{R}_+^d$, where $b$ and $\sigma$ satisfy non-negative
and suitable regularity conditions, and $X_t=(X_t^1,\cdots,X_t^d)$
represents $d$ populations. The branching rate of the $i$th
population of $X$ is a function $(\sigma_i)$ of the mass of  $d$
populations.

Infinite systems of mutually catalytic branching and cyclically
catalytic branching diffusions with $d\geq2$ and a linear
interaction between the components have been extensively studied in
Dawson and Perkins (1998), Dawson et al (2003) and Fleischmann and
Xiong (2001). Uniqueness for these systems follows from Mytnik's
self-duality; see Mytnik(1998). But this argument works only for
$d=2$. Swart (2004)  described  a new way to generalize mutually
catalytic branching diffusion to the case $d>2$, but the set-up
there was rather special. Moreover, in all of the work mentioned
above only  weak uniqueness has been obtained.

Dawson et al (2005) studied the \textit{strong uniqueness} problem
for cyclically catalytic branching diffusions in the simplified
space-less case. They addressed \textit{pathwise uniqueness} for the
SDE
 \beqn
 \label{0.2}
 dX_t^i=\alpha_iX_t^idt+\sqrt{\gamma_iX_t^{i-1}X_t^i}dB_t^i,~~~t>0,~~~i=1,\cdots,d,
 \eeqn
where $(\gamma_i,~i=1,\cdots,d)$ are strictly positive constants. In
this note, we study a slightly more general form of ($\ref{0.2}$).
Fix an integer $d\geq 1$, let $I^d:=\{1,2,\cdots,d\}$ and
$\textbf{f}:=(f_1,f_2,\cdots,f_d)$, where for each $i\in I^d$
,$$f_i:\mathbb{R}^d\mapsto [0,\infty)$$ is a  continuous function.
Consider the following stochastic differential equation
 \begin{equation}
 \label{1.1}
 \left\{\begin{array}{lll}
 dX_t^i&=&\alpha_iX_t^idt+\sqrt{f_i(X_t)X_t^i}dB_t^i,~~~t>0,~~~i\in
 I^d,\cr
 X_0&=&\textbf{a}=(a_1,a_2,\cdots,a_d)\in \mathbb{R}_+^d
 \end{array}\right.
 \end{equation}
for a diffusion process $\textbf{X}=(X^i)_{i\in I^d}$ in
$\mathbb{R}_+^d$. Here $(\alpha_i)_{i\in I^d}$ are real constants
and $\textbf{B}=(B_t^1,B_t^2\cdots,B_t^d)$ is a
$\mathbb{R}^d$-valued standard Brownian motion. Our main purpose is
to establish the pathwise uniqueness for  equation ($\ref{1.1}$).
The idea behind this uniqueness is as follows. Indeed, away from the
zero boundary, uniqueness holds by an ``extended Lipschitz
condition'' which was suggested by Fang and Zhang (2005). On the
other hand, once a component, say $X^k$, reaches zero, it is trapped
there. But after this trapping, the model simplifies drastically.
Then we can repeat the previous argument for the simplified model
and get the uniqueness result when the cycle is closed.

In section 2, we will describe the main results. The proof of the
uniqueness result will be given in section 3. With $C$ we denote a
positive constant which might change from line to line. For
$x\in\mathbb{R}^d$, let $|x|$ denote the Euclidean norm. For the
definitions of weak solution, strong solution, weak uniqueness,
pathwise uniqueness, explosion, etc., see Ikeda and Watanabe (1989)
for example.

\section{Main results}

\bigskip

\begin{theorem}
\label{SDE} Let $r$ be a strictly positive $C^1$-function
defined on an interval $(0,c_0]$, satisfying\\
 (i) $\liminf_{s\rightarrow 0}r(s)>0$,\\
 (ii) $\lim_{s\rightarrow0}\frac{sr'(s)}{r(s)}=0$, \\
 (iii) $\int_0^a\frac{ds}{sr(s)}=+\infty~ \textrm{for any}~a>0$.\\
Assume that equation ($\ref{1.1}$) has no explosion and for
$|x-y|\leq c_0$,
 \beqn
 \label{1.2}
 |\textbf{f}(x)-\textbf{f}(y)|^2\leq C|x-y|^2r(|x-y|^2).
 \eeqn
Then the pathwise uniqueness holds for stochastic differential
equation ($\ref{1.1}$) if one of the following conditions holds:

\begin{enumerate}
\item[(I)]
\textsl{$f_i(x)> 0$ for all $i\in I^d$ and $x\in \mathbb{R}_+^d$;}

\item[(II)]
\textsl{$\{x:f_i(x)=0\} \subset \partial \mathbb{R}_+^d$ for all
$i\in I^d$ } \textsl{and if $f_i(z)=0$ for some $z=(z_1,\cdots,z_d)
\in \partial \mathbb{R}_+^d$ satisfying $z_{i_1} =\cdots =
z_{i_k}=0$ and $z_j>0$ for $j\notin \{i_1,\cdots,i_k\}$, then
$f_i(x)=0$ for all $x= (x_1,\cdots,x_d) \in \partial\mathbb{R}_+^d$
satisfying $x_{i_1}=\cdots=x_{i_k}=0$.}
\end{enumerate}
\end{theorem}
\smallskip

\begin{remark} By letting $f_i(x)=\gamma_ix_i$, we shall see that
equation ($\ref{0.2}$) is a special case of equation ($\ref{1.1}$).
Therefore, the  result of uniqueness of Dawson et al (2005) is a
consequence of Theorem \ref{SDE}.
\end{remark}

\begin{remark} Functions $r(s)=\log 1/s,~r(s)=\log 1/s\cdot
\log\log 1/s,\cdots$ are typical examples satisfying the conditions
(i)-(iii) in Theorem \ref{SDE}.
\end{remark}
\begin{remark} Define the function V on $\mathbb{R}$ by the
series
$$
V(x)=\sum_{k=1}^{+\infty}\frac{|\sin kx|}{k^2}.
$$
For $X=(x_1,x_2)\in \mathbb{R}^2$ and  $\theta_1>0,~\theta_2>0$,
define $$f_i(X)=V(x_1)+V(x_2)+\theta_i,~~i=1,2.$$ According to
Example 2.8 of Fang and Zhang (2005), the bounded function
$\textbf{f}=(f_1,f_2)$ is of linear growth and satisfies the
inequality ($\ref{1.2}$) with $r(s)=\log 1/s$ and the condition (I)
in Theorem \ref{SDE}.
\end{remark}
\bigskip

\noindent Theorem \ref{SDE} is based on  the non-explosion
assumption of the equation ($\ref{1.1}$). The following proposition
gives a sufficient condition for non-explosion.

\bigskip

\begin{proposition}\label{prop}Let $\rho$ be a strictly
 positive $C^1$-function
defined on an interval $[K,+\infty)$ satisfying (i)
$\lim_{s\rightarrow\infty}\rho(s)=+\infty,$ (ii)
$\lim_{s\rightarrow\infty}\frac{s\rho'(s)}{\rho(s)}=0$ and
 (iii)$~\int_K^{+\infty}\frac{ds}{s\rho(s)+1}=+\infty.$\\
Assume that for $|x|\geq K$,
$$\sum_{i=1}^df_i^2(x)\leq C[|x|^2\rho(|x|^2)+1].$$
Then the equation ($\ref{1.1}$) has no explosion.
\end{proposition}
\smallskip

\begin{remark} Functions $\rho(s)=\log s,~\rho(s)=\log s\cdot
\log\log s,\cdots$ are typical examples satisfying the conditions
(i)-(iii) in Proposition \ref{prop}.
\end{remark}
\smallskip
\noindent\textbf{Proof}. By Theorem A of Fang and Zhang (2005), the
desired result is obvious. $\Box$

\bigskip

\section{Proof of Theorem \ref{SDE}}According to condition (i) on the
function $r$, we can assume that there exists a  constant $C_1>0$
such that $r(\zeta)\geq 1/C_1$ for all $\zeta\in(0,c_0]$. Let
$\dz>0$, we define
 $$
 \phi_{\dz}(\zeta)=\int_0^{\zeta}\frac{ds}{sr(s)+\dz}~~\textrm{and}~~
 \Phi_{\dz}(\zeta)=e^{\phi_{\dz}(\zeta)},~~\zeta\geq 0.
 $$
By condition (iii) on $r$, we see that $\Phi_0(\zeta)=+\infty$ for
all $\zeta>0$. We have
 $$
 \Phi_{\dz}'(\zeta)=\frac{\Phi_{\dz}(\zeta)}{\zeta r(\zeta)+\dz}
 ~~\textrm{and}~~\Phi_{\dz}''(\zeta)=
 \frac{1-r(\zeta)-\zeta r'(\zeta)}{(\zeta r(\zeta)+\dz)^2}\Phi_{\dz}(\zeta).
 $$
By  conditions (i) and (ii) on the function $r$, there exists a
constant $C_2>0$ such that
$$
|1-r(\zeta)-\zeta r'(\zeta)|\leq C_2r(\zeta).
$$
So that
 \beqn
 \label{1.3}
 \Phi_{\dz}''(\zeta)\leq C_2\frac{\Phi_{\dz}(\zeta)r(\zeta)}{(\zeta
 r(\zeta)+\dz)^2}.
 \eeqn

\noindent Without loss of generality, we may assume that $d\geq 2$.
Fix $i\in I^d$ and $\textbf{a}\in \mathbb{R}_+^d$. Clearly, from
It\^{o}'s formula, $t\mapsto e^{-\alpha_it}X_t^i$ is a non-negative
martingale, implying that the zero state is a trap for this
martingale. Hence,  $X^i$ is trapped at 0 once it reaches 0.

Suppose we have two solutions
$\textbf{X}=(x_t^1,x_t^2,\cdots,x_t^d)$ and
$\textbf{Y}=(y_t^1,y_t^2,\cdots,y_t^d)$ to ($\ref{1.1}$) with the
same $\textbf{B}$ and satisfying
$\textbf{X}_0=\textbf{a}=\textbf{Y}_0.$ Let
$$\eta_t^i:=\sqrt{f_i(\textbf{X}_t)x_t^i}-
\sqrt{f_i(\textbf{Y}_t)y_t^i}, ~~~~ \xi_t^i:=|x_t^i-y_t^i|^2,
~~\textrm{ for}~~i\in I^d,$$ and
$$ \zeta_t:=|\textbf{X}_t-\textbf{Y}_t|^2.$$

\smallskip

\noindent Fix $\textbf{a}\in \mbb R_+^d$ and  $a_i>0$ for $i\in
I_d$. Let $\ez>0$ be such that $\ez<a_i<\ez^{-1},~i\in I^d$.
Introduce two stopping times
 \begin{eqnarray}
 \label{2.1}
 \tau_{\ez}^d:=\inf\bigg{\{}t>0: \exists  i\in
 I^d &\textrm{with}&
 f_i(\textbf{X}_t)\wedge f_i(\textbf{Y}_t)\wedge x_t^i\wedge
 y_t^i\leq \ez\cr
 &\textrm{or}&f_i(\textbf{X}_t)\vee f_i(\textbf{Y}_t)\vee
 x_t^i\vee y_t^i\geq \ez^{-1}\bigg{\}},
 \end{eqnarray}
and
$$
\tau:=\inf\{t>0:\zeta_t\geq c_0^2\}.
$$
\\
\begin{lemma}\label{le3} For any fixed $T>0$, we have $\textbf{X}=\textbf{Y}$ on $[0,
T\wedge \tau\wedge\tau_{\ez}^d]$.\end{lemma}

\smallskip

\noindent\textbf{Proof}. From equation ($\ref{1.1}$), we have
\begin{eqnarray}
\label{2.2} d\xi_t^i=2\alpha _i\xi_t^idt+
2(x_t^i-y_t^i)\eta_t^idB_t^i+(\eta_t^i)^2dt
\end{eqnarray}
and
 \beqn
 \label{2.3}
  d\la\xi^i,\xi^i\ra_t=4\xi_t^i(\eta_t^i)^2dt.
 \eeqn
According to ($\ref{1.2}$), for $s\leq \tau\we\tau_{\ez}^d$,
\begin{eqnarray}
\label{2.5} (\eta_s^i)^2\leq\frac{1}{2\ez^4}
  (|f_i(\textbf{X}_s)-f_i(\textbf{Y})_s|^2+\zeta_s)
  \leq\frac{1}{2\ez^4}(C\zeta_s r(\zeta_s)+\zeta_s),
\end{eqnarray}
where the first inequality is due to the elementary inequalities
$$
|\sqrt{bc}-\sqrt{de}|\leq \frac{1}{2\ez^2}(|b-d|+|c-e|)
~~\textrm{if}~~\ez\leq b,c,d,e\leq\ez^{-1}
$$
and $(g+h)^2\leq 2(g^2+h^2).$\\
Applying It\^{o}'s formula and according to ($\ref{2.2}$) and
($\ref{2.3}$), we have
\begin{eqnarray}
\label{2.4} \Phi_{\dz}(\zeta_{t\we\tau\we\tau_{\ez}^d})
&=&\Phi_{\dz}(\zeta_0)+\sum_{i=1}^d\int_0^{t\we\tau\we\tau_{\ez}^d}
   2\alpha_i\Phi_{\dz}'(\zeta_s)\xi_s^ids\cr
&&+\sum_{i=1}^d\int_0^{t\we\tau\we\tau_{\ez}^d}
  2\Phi_{\dz}'(\zeta_s)(x_s^i-y_s^i)\eta_s^idB_s^i
+\sum_{i=1}^d\int_0^{t\we\tau\we\tau_{\ez}^d}
  \Phi_{\dz}'(\zeta_s)(\eta_s^i)^2ds\cr
&&+\sum_{i=1}^d\int_0^{t\we\tau\we\tau_{\ez}^d}
2\Phi_{\dz}''(\zeta_s)\xi_s^i(\eta_s^i)^2ds\cr
&=&\Phi_{\dz}(\zeta_0)+I_1(t)+I_2(t)+I_3(t)+I_4(t)
\end{eqnarray}
respectively. For any $s\leq \tau\we\tau_{\ez}^d$,  by condition (i)
on $r$,
 \beqn
\sum_{i=1}^d2\alpha_i\Phi_{\dz}'(\zeta_s)\xi_s^i\leq
\frac{\alpha\zeta_s\Phi_{\dz}(\zeta_s)}{\zeta_sr(\zeta_s)+\dz}\leq
\alpha C_1\Phi_{\dz}(\zeta_s),
 \eeqn
where $\alpha=2\max_{1\leq i\leq d}|\alpha_i|$. According to
($\ref{2.5}$), we have
 \begin{eqnarray*}
 |\Phi_{\dz}'(\zeta_s)(x_s^i-y_s^i)\eta_s^i|^2
\leq\frac{1}{2\ez^4}\cdot\frac{\Phi_{\dz}^2(\zeta_s)(C\zeta_s^2r(\zeta_s)+\zeta_s^2)}{(\zeta_s
r(\zeta_s)+\dz)^2}\leq
\frac{(CC_1+C_1^2)}{2\ez^4}\sup_{0\leq\zeta\leq
c_0}\Phi_{\dz}(\zeta)^2<\infty
 \end{eqnarray*}
and
 \begin{eqnarray*}
 \Phi_{\dz}'(\zeta_s)(\eta_s^i)^2\leq
 \frac{\Phi_{\dz}(\zeta_s)}{2\ez^4}\cdot\frac{C\zeta_sr(\zeta_s)+\zeta_s}
 {\zeta_sr(\zeta_s)+\dz}\leq \frac{(C+C_1)}{2\ez^4}\Phi_\dz(\zeta_s).
 \end{eqnarray*}
On the other hand, by ($\ref{1.3}$) and ($\ref{2.5}$),
\begin{eqnarray*}
2\Phi_{\dz}''(\zeta_s)\xi_s^i(\eta_s^i)^2\leq
\frac{C_2}{\ez^4}\cdot\frac{\Phi_{\dz}(C\zeta_s^2r(\zeta_s)+\zeta_s^2)}{(\zeta_sr(\zeta)+\dz)^2}
\leq \frac{(CC_1C_2+C_1^2C_2)}{\ez^4}\Phi_{\dz}(\zeta_s).
\end{eqnarray*}
Therefore, $I_2(t)$ is a martingale and $\mathbb{E}(I_2(t))=0$. Let
$$
K:=\alpha C_1+\frac{d(C+C_1)}{2\ez^4}+\frac{dC_1C_2(C+C_1)}{\ez^4}.
$$
We have
$$
\mathbb{E}\bigg{(}\Phi_{\dz}(\zeta_{t\we\tau\we\tau_{\ez}^d})\bigg{)}\leq
\Phi_{\dz}(\zeta_0)+K\int_0^t
\mathbb{E}\bigg{(}\Phi_{\dz}(\zeta_{s\we\tau\we\tau_{\ez}^d})\bigg{)}ds.
$$
Thanks to Gronwall's inequality, we have that, for all $t>0$,
$$
\mathbb{E}\bigg{(}\Phi_{\dz}(\zeta_{t\we\tau\we\tau_{\ez}^d})\bigg{)}\leq
\Phi_{\dz}(\zeta_0)e^{Kt}.
$$
Letting $\dz \downarrow0$ in the above inequality, we find that
$$
\mathbb{E}\bigg{(}\Phi_0(\zeta_{t\we\tau\we\tau_{\ez}^d})\bigg{)}\leq
e^{Kt}.
$$
By the continuity of the samples and the fact that
$\Phi_0(\zeta)=+\infty$ for $\zeta>0$, we can get almost surely
$$
\zeta_{t\we\tau\we\tau_{\ez}^d}=0~~\textrm{for all }t>0.
$$
This yields the desired result. $\square$

\medskip

As $\ez\downarrow0$, we have the non-decreasing convergence of
$\tau_\ez^d$ to some stopping time $\tau^d\leq \infty$. On
$\{\tau^d=\infty\}$, we clearly have $\zeta_{t\we\tau}=0$, and by
the continuity of samples and the definition of $\tau$, we have that
almost surely  for all  $t>0$, $\zeta_t=0$. Thus, we may assume that
$\mathbb{P}(\tau^d<\infty)>0$. And on $\{\tau^d<\infty\}$, by
similar reasoning, $\zeta_{t\we\tau\we\tau^d}=0$ implies
 \beqn
 \label{2.6}
 \zeta_{t\we\tau^d}=0.
 \eeqn
Now, we are in position to complete the proof of Theorem \ref{SDE}.

Let $\Omega_0^d=\{\tau^d<\infty\}$. We only need to show that for
all $t\geq0$  $$\zeta_t=0, \textrm{ on }\Omega_0^d.$$

\noindent\textit{Case I}: Since the stochastic differential equation
(\ref{1.1}) has no explosion, for each $\omega\in\Omega_0^d$, there
exists a $k\in I^d$ such that
$x^k_{\tau^d}(\omega)=0=y^k_{\tau^d}(\omega)$. For $k\in I^d$,
define
$$
\Omega_k^d=\{\omega\in\Omega_0^d~:~ x^k_{\tau^d}=0=y^k_{\tau^d}\}.
$$
We shall see from the following argument that there is no loss of
generality if we
 assume that for each $k\in I^d$, ${\mbb P}(\Omega_k^d)>0$ and
 $\Omega_i^d\cap\Omega_j^d=\emptyset$ for $i\neq j$. Note that
 $$\Omega_0^d=\cup_{k=1}^d\Omega_k^d.$$
  Since 0 is a trap, for $\omega\in\Omega_k^d$,
 $x_t^k(\omega)=0=y_t^k(\omega)$, for any $t>\tau^d(\omega)$.
 Together with ($\ref{2.6}$),
we have for any $t\geq 0$, $x_t^k=y_t^k$ on $\Omega_k^d$. To show
that $\bf X$ is pathwise uniquely determined, we only need to show
that for each $l\in I^d\backslash\{k\}$,
$$\xi_{\cdot}^l=0,\textrm{ on } \Omega_k^d.$$
 For $\ez>0$, define two new stopping times such that on
$\Omega_k^d$
\begin{eqnarray}
 \label{2.7}
 \tau_{\ez}^{d-1}:=\inf\bigg{\{}t>0: \exists~  i\in
 I^{d}\backslash\{k\} &\textrm{with}&
 f_i(\textbf{X}_t)\wedge f_i(\textbf{Y}_t)\wedge x_t^i\wedge
 y_t^i\leq \ez\cr
 &\textrm{or}&f_i(\textbf{X}_t)\vee f_i(\textbf{Y}_t)\vee
 x_t^i\vee y_t^i\geq \ez^{-1}\bigg{\}}
 \end{eqnarray}
 and
 \begin{eqnarray}
 \label{2.8}
 \tau^{\ez}_{1}:=\inf\bigg{\{}t>0: ~x_t^k\geq\ez^{-1}&\textrm{or}&
 f_k(\textbf{X}_t)\wedge f_k(\textbf{Y}_t)\leq \ez\cr
 &\textrm{or}&f_k(\textbf{X}_t)\vee f_k(\textbf{Y}_t)\geq \ez^{-1}\bigg{\}}
 \end{eqnarray}
and $\tau_1^{\ez}=\tau_{\ez}^{d-1}=\tau_{\ez}^d$ on
$\{\tau^d=\infty\}$. Note that $\tau_1^{\ez}\rightarrow\infty$ as
$\ez\rightarrow0$. Let
$\gamma^{\ez}=\tau_1^{\ez}\we\tau_{\ez}^{d-1}$. For
$s\leq\tau\we\gamma^{\ez}$,
 $$
 (\eta_s^i)^2\leq\frac{1}{2\ez^4}
  (|f_i(\textbf{X}_s)-f_i(\textbf{Y})_s|^2+\zeta_s)
  \leq\frac{1}{2\ez^4}(C\zeta_s r(\zeta_s)+\zeta_s), \textrm{ on }
  \Omega\backslash\Omega_i^d$$
and
$$
  (\eta_s^i)^2=(x_s^i)^2(\sqrt {f_i(\textbf{X}_s)}-\sqrt{f_i(\textbf{Y})_s})^2
  \leq\frac{1}{4\ez^6}|f_i(\textbf{X}_s)-f_i(\textbf{Y})_s|^2
  \leq\frac{1}{4\ez^6}C \zeta_s r(\zeta_s), \textrm{ on } \Omega_i^d.
 $$
Thus
 \beqn
 \label{3.15}
 (\eta_s^i)^2\leq
  \left(\frac{1}{2\ez^4}\vee\frac{1}{4\ez^6}\right)
  (C\zeta_s r(\zeta_s)+\zeta_s).
 \eeqn
By the same argument as in Lemma \ref{le3}, we have that almost
surely
$$
\zeta_{t\we\tau\we\gamma^{\ez}}=0,\textrm{ for~ all } t\geq0.
$$
And, we also have that there exists a stopping time
$\tau^{d-1}\leq\infty$ such that
$\tau_{\ez}^{d-1}\rightarrow\tau^{d-1}$. Note that
$\{\tau^d=\infty\}\subset\{\tau^{d-1}=\infty\}$ and  we have for
$t\geq0$,
$$
\zeta_t=0, \textrm{ on } \{\tau^{d-1}=\infty\}.
$$
Let $\Omega_0^{d-1}=\{\tau^{d-1}<\infty\}$. For each
$\omega\in\Omega_0^{d-1}\cap\Omega_k^d$, there exists a $l\in
I^d\backslash\{k\}$ such that
$x^l_{\tau^{d-1}}(\omega)=0=y^l_{\tau^{d-1}}(\omega)$. For $l\in
I^d\backslash\{k\}$, define
$$
\Omega_{kl}^{d-1}:=\{\omega\in\Omega_0^{d-1}\cap\Omega_k^d:~
x^l_{\tau^{d-1}}(\omega)=0=y^l_{\tau^{d-1}}(\omega)\}.
$$
Note that
$$
\bigcup_{k\in I^d}\bigcup_{l\in I^d\backslash
\{k\}}\Omega_{kl}^{d-1}=\{\tau^{d-1}<\infty\}.
$$
Since $0$ is a trap, we see for all $t\geq0$,
$$
x_t^k=y_t^k \textrm{ and }x_t^l=y_t^l,\textrm{ on
}\Omega_{kl}^{d-1}.
$$
Note that if $d=2$, then we are done. Next, we still assume that for
$i\neq j$, $\Omega_{ki}^{d-1}\cap\Omega_{kj}^{d-1}=\emptyset$. For
$\ez>0$, define two new stopping times such that on
$\Omega_{kl}^{d-1}$
 \begin{eqnarray}
 \label{3.7}
 \tau_{\ez}^{d-2}:=\inf\bigg{\{}t>0: \exists~  i\in
 I^{d}\backslash\{k,l\} &\textrm{with}&
 f_i(\textbf{X}_t)\wedge f_i(\textbf{Y}_t)\wedge x_t^i\wedge
 y_t^i\leq \ez\cr
 &\textrm{or}&f_i(\textbf{X}_t)\vee f_i(\textbf{Y}_t)\vee
 x_t^i\vee y_t^i\geq \ez^{-1}\bigg{\}}
 \end{eqnarray}
 and
 \begin{eqnarray}
 \label{3.8}
 \tau^{\ez}_{2}:=\inf\bigg{\{}t>0: \exists~ i\in\{k,l\} \textrm{ with }~x_t^i\geq\ez^{-1}&\textrm{or}&
 f_i(\textbf{X}_t)\wedge f_i(\textbf{Y}_t)\leq \ez\cr
 &\textrm{or}&f_i(\textbf{X}_t)\vee f_i(\textbf{Y}_t)\geq \ez^{-1}\bigg{\}}
 \end{eqnarray}
and $\tau_2^{\ez}=\tau_{\ez}^{d-2}=\tau_{\ez}^{d-1}$ on
$\{\tau^{d-1}=\infty\}$. Repeat the previous argument for
$\tau_{\ez}^{d-2}$ instead of $\tau_{\ez}^{d-1}$ and $\tau^{\ez}_2$
instead of $\tau^{\ez}_1$. Then we get a new partition  on $\Omega$,
say $\{\Omega_i\}_{i=0}^n$, and  for each $\Omega_i$ there exist at
least three components of $\textbf X$ such that they are pathwise
uniquely determined on $\Omega_i$. By this way the argument can be
repeated until the cycle is closed. We conclude that pathwise
uniqueness holds for equation ($\ref{1.1}$).

\textit{Case II}: By similar reasoning, for $\omega\in\Omega_0^d$,
we have that there exists a $k\in I^d$ such that
$x_{\tau^d}^k(\omega)=0=y_{\tau^d}^k(\omega)$ and a $l\in I^d$ such
that
$f_l(\textbf{X}_{\tau^d}(\omega))=0=f_l(\textbf{Y}_{\tau^d}(\omega))$
[$k$ may be equal to $l$]. Using the previous argument, we find that
$$x_t^k=y_t^k,~~~~t\geq 0.$$
Since 0 is a trap, and by the condition on $f_l$, 0 is also a trap
for random processes $\{f_l(\textbf{X}_t): t\geq 0\}$ and
$\{f_l(\textbf{Y}_t): t\geq 0\}$. This implies, after the trapping
event, $dx_t^l=\alpha_lx_t^ldt$ and $dy_t^l=\alpha_ly_t^ldt$. That
is
$$x_t^l=x_{\tau^d}^le^{\alpha_lt}=y_{\tau^d}^le^{\alpha_lt}=y_t^l,~~~~t\geq \tau^d.$$
Also, for all $t\geq0$, $f_l({\bf X}_t)=f_l({\bf Y}_t)$ on
$\Omega_{lf}^d:=\{\omega\in\Omega_0^d:~f_l(\textbf{X}_{\tau^d})=0=f_l(\textbf{Y}_{\tau^d})\}$.
Note that $\{x:f_l(x)=0\} \subset \partial \mathbb{R}_+^d$. Define
$\Omega_k^d$ as that in previous case. Assume that
$\Omega_i^d\cap\Omega_j^d=\emptyset$ and
$\Omega_{if}^d\cap\Omega_{jf}^d=\emptyset$ for $i\neq j$. For
$\ez>0$, introduce two stopping times such that on
$\Omega_k^d\cap\Omega_{lf}^d$
\begin{eqnarray*}
 \tau_{\ez}^{k,lf}:=\inf\bigg{\{}t>0: \exists ~ i\in
 I^{d}\setminus\{k,l\} &\textrm{with}&
 f_i(\textbf{X}_t)\wedge f_i(\textbf{Y}_t)\wedge x_t^i\wedge
 y_t^i\leq \ez\cr
 &\textrm{or}&f_i(\textbf{X}_t)\vee f_i(\textbf{Y}_t)\vee
 x_t^i\vee y_t^i\geq \ez^{-1}\bigg{\}}
 \end{eqnarray*}
 and
\begin{eqnarray*}
  \tau^{\ez}_{k,lf}:=\inf\bigg{\{}t>0: ~x_t^k\vee x_t^l\vee
   f_k({\bf X}_t)\vee f_l({\bf X}_t)\geq\ez^{-1} \bigg{\}}
 \end{eqnarray*}
 and on $\Omega_k^d\cap(\cup_{l\in I^d}\Omega_{lf}^d)^c$, they are
 defined by ($\ref{2.7}$) and ($\ref{2.8}$) respectively and they equal to
 $\tau_{\ez}^{d}$ on
$\{\tau^{d}=\infty\}$. Then the argument would be exactly parallel
to that used in Case I. We omit it here and get the pathwise
uniqueness of $\textbf{X}$. This completes the proof of the theorem.
$\square$

\end{document}